\newcommand{\N}{\mathbb{N}}               

\newcommand{\R}{\mathbb{R}}
\mathchardef\varepsilon="010F
\mathchardef\epsilon="0122
\mathchardef\vartheta="0112
\mathchardef\varrho="011A\mathchardef\rho="0125
\mathchardef\varphi="011E     
\mathchardef\phi="0127
\renewcommand \emptyset \varnothing
\documentclass[10pt]{article}            
\usepackage[all]{xy}
\usepackage[latin1]{inputenc}        
\usepackage[dvips]{graphics}
\usepackage[dvips]{graphicx}
\usepackage{amsfonts}
\usepackage{amssymb}
\usepackage{amsmath}
\usepackage{amstext}
\usepackage{amsbsy}
\usepackage{amsopn}
\usepackage{amsthm}
\usepackage{amscd}
\usepackage{amsxtra}
\usepackage{mathrsfs}
\usepackage{upref}
\date\today                     
\title{On some nonlinear partial differential equations involving the 1-Laplacian}
\author{\sc Mouna Kra\"iem \\ 
{\small Université de Cergy Pontoise , Departement de Mathematiques,}\\ 
{\small 2, avenue Adolphe Chauvin, 95302 Cergy Pontoise Cedex, France} \\ 
{\small e-mail: Mouna.Kraiem@math.u-cergy.fr}}
\numberwithin{equation}{section}                
\theoremstyle{plain} 
\newtheorem{theorem*}{\bf Theorem}
\newtheorem{theorem}{\bf Theorem}
\newtheorem{lemma}{\bf Lemma}
\newtheorem{proposition}{\bf Proposition}
\newtheorem*{proposition*}{\bf Proposition}

\newtheorem*{corollary*}{\bf Corollary}
\theoremstyle{definition}
\newtheorem*{definition*}{\bf Definition}
\newtheorem*{definitions*}{\bf Definitions}

\newtheorem*{example*}{\bf Example}
\theoremstyle{remark}
\newtheorem{remark}{\bf Remark}

\begin{document}  
\maketitle

\begin{center}
{\bf Abstract}
\end{center}
Let $\Omega $ be a smooth bounded domain in $\R^N, N>1$ and let $n\in \N^*$. 
We prove here the existence of nonnegative solutions $u_n$ in $BV(\Omega)$, to the problem
$$(P_n) \begin{cases}-{\rm div} \sigma +2n \left(\int_ \Omega u -1 \right) \ {\rm sign}^+ \ (u)=0\ \quad \text{in} \ \Omega,\\\sigma \cdot \nabla u= |\nabla u|  \quad \text{in} \ \Omega, \\u \  \text{\rm is not identically zero}, -\sigma \cdot \overrightarrow {n} u=u  \quad \text{on}\ \partial\Omega,\end{cases}$$
where $\overrightarrow {n}$ denotes the unit outer normal to $\partial\Omega$, and ${\rm sign}^+(u)$ denotes some $L^{\infty}(\Omega)$ function defined as:
$${\rm sign}^+ \ (u). u =u^+, \  0 \leq {\rm sign}^+(u) \leq 1.$$ 
Moreover, we prove the tight convergence of $u_n$ towards one of the first eingenfunctions for the first $1-$Laplacian Operator $-\Delta_1$ on $\Omega$ when $n$ goes to $+\infty$. 

\begin{center}
{\bf Résumé}
\end{center}
Soit $\Omega $ un domaine borné et lisse dans $\R^N, N>1$ et soit $n\in \N^*$.
On montre dans ce papier l'existence de solutions positives $u_n$ dans  $BV(\Omega)$, au probléme 
$$(P_n) \begin{cases}-{\rm div} \sigma +2n \left(\int_ \Omega u -1 \right) \ {\rm sign}^+ \ (u)=0\ \quad \text{in} \ \Omega,\\\sigma \cdot \nabla u= |\nabla u|  \quad \text{in} \ \Omega, \\u \  \text{\rm is not identically zero}, -\sigma \cdot \overrightarrow {n} u=u  \quad \text{on}\ \partial\Omega,
\end{cases}$$
où $\overrightarrow {n}$ est le vecteur normal sortant de $\partial\Omega$, et 
${\rm sign}^+(u)$ est une fonction dans $L^{\infty}(\Omega)$ definie par:

$${\rm sign}^+ \ (u). u =u^+, \  0 \leq {\rm sign}^+(u) \leq 1.$$ 
De plus, on montre la convergence de $u_n$ vers une des valeurs premiéres de l'opérateur 
$1-$Laplacian  $-\Delta_1$ sur  $\Omega$ quand  $n$ tend vers $+\infty$. \\
{\small {\bf  Mots clés:}} Fonctions BV, 1-Laplacien.

\section{Introduction}

Recent works about the operator $1$-Laplacian revealed the existence of a least eigenvalue for this  operator on a bounded smooth set. 
More precisally this first eigenvalue is well defined as the infimum of 
$$\inf_{u\in W_0^{1,1}(\Omega), \ \int_\Omega |u|=1} \int_\Omega |\nabla u|.$$
This value denoted as $\lambda_1$ is positive by Poincaré 's inequality.
Unfortunately, since $W^{1,1}(\Omega)$ is not a reflexif space, it is not possible to prove the existence of solutions in $W_0^{1,1}(\Omega)$.
The convenient space in which one must look for a minimizer is the space $BV(\Omega)$ which is the weak closure of $W^{1,1}(\Omega)$. 
Moreover, since the trace map which is well defined on $BV(\Omega)$  is not weakly continuous, one is lead to replace the problem by the relaxed following form
$$\inf_{u\in BV(\Omega), \int_\Omega |u|=1} \int|\nabla u|+\int_{\partial\Omega} |u|. $$
This problem has an infimum equal to $\lambda_1$. 
Classical arguments in the theory of $BV$ functions allow then to prove the existence of a minimizer. 
Moreover using either the fact that there exist non negative solutions and duality in convex analysis, or using an approximation with the more regular problem 
$$\inf_{u\in W_0^{1,1+\epsilon} (\Omega ), \int |u|=1} \int_\Omega|\nabla u|^{1+\epsilon},$$

one can obtain that the solution $u$ satisfies the singular PDE 
$$-{\rm div} (\frac{\nabla u}{|\nabla u|}) = \lambda_1,$$
in a sense which must of course be precised, and is detailed in the present paper. 
Let us note that it is proved in \cite{DF3} that there are caracteristic functions of sets which are solutions, they are consequently called eigensets. 
Another approach is used in \cite{ACCh1} \cite{ACCh2}, where the authors use the concept of Cheeger sets. In these papers, the authors present a remarquable construction of eigenset in the case $N=2$ and for convex sets $\Omega$. 
Among their results there is  the uniqueness of eigen sets in this case. Our aim in the present paper is to propose an approach of the first eigenvalue and first eigenfunction, using a penalization method. 
This method has an obvious numerical advantage : The constraint $\int_{\Omega}|u|=1$ has a higher coast than the introduction of the penalization term $n(\int_{\Omega}|u|-1)^2$.
In the same time one gets a new proof of the existence of nonnegative function. 

\section{Preliminaries}

Let $\Omega $ be a smooth bounded domain in  $\R^N, N>1$, whose boundary is piecewise $\mathcal C^1$, and let us define for all $n\in \N^*$,  the following functional
$$I_{n,0}(u)=  \int_\Omega |\nabla u| +  n \left(\int_\Omega|u| -1 \right)^2 \ .$$  
It is clear, using Poincaré's inequality, that there exists some constant $c>0$, such that for all $u \in W^{1,1}_0(\Omega)$ and for all $n\in \N^*$,

\begin{equation}\label{L1.1}
\int_\Omega |\nabla u| +  n \left( \int_\Omega|u| -1\right)^2 \ \geq c \ \|u\|_{ W^{1,1}(\Omega)}.
\end{equation}
We look for $u \in W^{1,1}_0(\Omega)$, nonnegative  which satisfies:

\begin{equation}\label{L1.2}
\begin{cases}-{\rm div} \sigma +2n \left( \int_ \Omega u -1 \right) \ {\rm sign}^+ \ (u)=0 \quad \text{in} \ \Omega,\\
\sigma \in L^{\infty}( \Omega,\R^n),\\
\sigma \cdot\nabla u= |\nabla u|  \quad \text{in} \ \Omega, \\
u \  \text{\rm is not identically zero} , u=0 \quad \text{on}\ \partial\Omega,
\end{cases}
\end{equation}
In order to find solutions to (\ref{L1.2}) one can consider the following minimisation problem

\begin{equation}\label{L1.3}
\inf_{\begin{smallmatrix} 
u \in W^{1,1}_0(\Omega)  
\end{smallmatrix}} 
\left\{ \int_\Omega |\nabla u| +  n \left( \int_\Omega|u| -1 \right)^2 \right\},  
\end{equation}
We denote by $\lambda_n(\Omega)$ the value of this infimum. In the following, we shall prove that, if $u_n \in W^{1,1}_0(\Omega)$ realizes the minimum defined in (\ref{L1.3}), it is a non trivial solution of (\ref{L1.2}).\\
Since classical methods in the calculus of variations cannot be applied to solve (\ref{L1.3}), we approximate it by the following formulation: for $\epsilon>0$, we  define

$$\lambda_{n,\epsilon}(\Omega)=
\inf_{\begin{smallmatrix} 
u \in W^{1,1+\epsilon}_0(\Omega)
\end{smallmatrix}} 
\left\{\int_\Omega |\nabla u |^{1+\epsilon} +  n \left(\int_\Omega|u|^{1+\epsilon} -1 \right)^2 \right \}.$$ 
Note that a non trivial, nonnegative minimizer $u_{n,\epsilon}$ for this problem solves the following partial differential equation:

$$\begin{cases}
-{\rm div} (|\nabla u_{n,\epsilon}|^{\epsilon-1} \nabla u_{n,\epsilon}) +2n \left(\int_ \Omega u_{n,\epsilon}^{1+\epsilon} -1 \right) u_{n,\epsilon}^\epsilon =0  \quad \text{in} \ \Omega, \\u_{n,\epsilon} \in W^{1,1+\epsilon}_0(\Omega).  
\end{cases}$$
A solution of (\ref{L1.2}) will be obtained, letting $\epsilon$ tend to $0$.
Let us observe that, in particular, regularity and other properties of $u_n$ can be derived from a priori estimates on $u_{n,\epsilon}$.
Of course, passing to the limit when $\epsilon \rightarrow 0$ will lead us to 
consider $BV(\Omega)$ in place of $W^{1,1}(\Omega)$, 
and to give sense to some expressions as $\sigma_n.\nabla u_n$ when $\nabla u_n$ 
is only a measure, and $\sigma_n \in L^{\infty}(\Omega)$, 
${\rm div}\sigma_n \in L^N(\Omega)$.

As the ``limit'' will be obtained by weak convergence in  $BV(\Omega)$, 
we shall be led to overcome the lack of continuity of the trace map of this space for weak topology. This can be done by introducing the concept of ``relaxed problem'': these problems are used in the theory of minimal surfaces and plasticity, and with a slightly different meaning, in the theory of weakly harmonic functions. 
Here the relaxed problem is defined as: 

\begin{equation}\label{L1.4}
\inf_{\begin{smallmatrix} 
u \in BV(\Omega) 
\end{smallmatrix}} 
\left\{\int_\Omega |\nabla u| +\int_{\partial \Omega} |u|+  n \left( \int_\Omega|u| -1\right)^2 \right\}.
\end{equation}
We shall prove in the sequel that (\ref{L1.4}) has the same infimum as (\ref{L1.3}) and that it possesses a solution $u$ in $BV(\Omega)$ which satisfies an equation as (\ref{L1.2}), extended to $BV$-functions, as it is done in \cite{DF6}. Of course, to prove this, and as we pointed out before, one must give sense to the product $''\sigma.\nabla u''$ when $\nabla u$ is only a measure.

\begin{proposition}\label{prop:L1}
Let  $\sigma$ be  in $ L^{\infty}(\Omega)$ and ${\rm div}\sigma \in L^N(\Omega)$, and define the distribution $\sigma.\nabla u$ by the formula : For  $\phi \in \mathcal D(\Omega)$, $u \in BV(\Omega)$,

\begin{equation}\label{L1.5}
\langle \sigma \cdot \nabla u, \phi \rangle=-\int_\Omega {\rm div}\sigma \  u \ \phi-\int_\Omega \sigma \cdot \nabla \phi \ u .
\end{equation}
Then

$$| \langle \sigma \cdot \nabla u, \phi \rangle| \leq | \sigma|_{\infty} \langle| \nabla u|,| \phi| \rangle .$$
In particular, $ \sigma \cdot \nabla u$ is a bounded measure on $\Omega$ which is absolutely continuous with respect to $|\nabla u|$.
In addition, if $\phi \in \mathcal C(\overline{\Omega}) \cap  C^1(\Omega)$, the following Green's Formula holds

\begin{equation}\label{L1.6}
\langle \sigma \cdot \nabla u, \phi \rangle=-\int_\Omega {\rm div}\sigma \  u \ \phi-\int_\Omega \sigma \cdot  \nabla \phi \ u +\int_{\partial \Omega} \sigma \cdot \overrightarrow {n} u \ \phi.
\end{equation}
Suppose that $U \in BV(\R^N- \overline{\Omega})$, and define for $u \in BV(\Omega)$ the functional $\tilde u$ as:

$$\tilde u =
\begin{cases}
u \ {\rm in} \ \Omega,\\
U \ {\rm in} \ \R^N- \overline{\Omega},\\
\end{cases}$$
then $\tilde u \in BV(\R^N)$ and
$$\nabla \tilde u = \nabla u \chi_{\Omega} +\nabla U \chi_{(\R^N- \overline{\Omega})}+(U-u)_ {| \Omega} \ \delta_{\partial \Omega},$$
where $U_ {| \Omega}$ and $u_ {| \Omega}$ denote the trace of $U$ and $u$ on $\partial \Omega$, $\delta_{\partial \Omega}$ denotes the uniform Dirac measure on  $\partial \Omega$ and $\overrightarrow {n}$ is the unit outer normal to $\partial \Omega$. Finally,  we introduce the measure $\sigma \cdot \nabla\tilde u$ on $\overline{\Omega}$ by the formula

$$(\sigma \cdot \nabla \tilde u)=(\sigma \cdot \nabla u)\chi_{\Omega} +\sigma \cdot\overrightarrow {n}(U-u) \ \delta_{\partial \Omega} .$$
Then $\sigma \cdot \nabla \tilde u$ is absolutely continuous with respect to $|\nabla \tilde u|$, with the inequality$$|\sigma \cdot \nabla \tilde u| \leq |\sigma|_{\infty} |\nabla  \tilde u|.$$
\end{proposition}
\noindent For a proof the reader can consult (\cite{DF2},\cite{Ko.Tem},\cite{Str.Tem}).

\section{Existence's Theorem}

\begin{theorem}\label{th:L1}
Let $\Omega$ be a bounded domain in $\R^N, N>1$, which is a piecewise $\mathcal C^1$set.
Let $\lambda_n$ be defined as in {\rm (\ref{L1.3})}. There exists a nonnegative solution $u_n \in BV(\Omega)$ to the problem  {\rm(\ref{L1.2})}which realizes the following partial differential equation: 

\begin{equation}\label{L2.1}
\begin{cases}
-{\rm div} \sigma_n +2n \left(\int_ \Omega u_n -1\right) {\rm sign}^+ (u_n) =0\quad \text{\rm in} \ \Omega, \\\sigma_n \in  L^{\infty}( \Omega,\R^N), |\sigma_n|_{\infty} \leq 1,\\\sigma_n \cdot \nabla u_n= |\nabla u_n|  \quad \text{\rm in} \ \Omega, \\
u \  \text{\rm is not identically zero}, 
-\sigma_n \cdot \overrightarrow {n} (u_n)= u_n 
\quad \text{\rm on}\ \partial\Omega,
\end{cases}
\end{equation}
$\overrightarrow {n}$ denotes the unit outer normal to $\partial\Omega$, 
and $\sigma_n \cdot\nabla u_n$ is the measure defined in proposition 
\ref{prop:L1}.
\end{theorem}

\begin{remark}
From proposition \ref{prop:L1} (with $U=0$) the conditions:

$$ \sigma_n \cdot \nabla u_n= |\nabla u_n|  \quad \text{\rm in} \ \Omega \ 
, \  \quad \ 
-\sigma_n \cdot \overrightarrow {n} (u_n)= u_n 
\quad \text{\rm on}\ \partial\Omega \ $$
are equivalent to
$$\sigma \cdot \nabla \tilde u =|\nabla \tilde u| \quad \text{\rm on}\ 
\overline{\Omega}$$
\end{remark}

\begin{remark}
Eq (\ref{L2.1}) can be written:

$${\rm div} \sigma_n= -2n \left(\int_ \Omega u_n -1\right) {\rm sign}^+ \ (u_n)$$
and since the right-hand side is an element of  $L^{\infty}(\Omega)$, 
$\sigma_n \cdot \nabla u_n$ is well-defined.
\end{remark}

\begin{proof}[Proof of Theorem \ref{th:L1}]

We split the proof of Theorem \ref{th:L1} into several steps:\\
{\bf Step 1}:\\
We begin by approximating (\ref{L1.4}) with the following minimization problem

\begin{equation}\label{L2.2}
\lambda_{n,\epsilon}(\Omega)=\inf_{\begin{smallmatrix} u \in 
W^{1,1+\epsilon}_0(\Omega)
\end{smallmatrix}} 
\left\{\int_\Omega |\nabla u |^{1+\epsilon} +  
n \left(\int_\Omega|u|^{1+\epsilon} -1 \right)^2 \right \},
\end{equation} 
where $\epsilon$ is some positive parameter. 
This problem can be solved by classical methods in the calculus of variations, 
since the compact embedding of $W^{1,1+\epsilon}_0(\Omega)$ into 
$L^q(\Omega)$ holds for all $q< (1+\epsilon)^*$.
Furthermore, there exists a nonnegative solution to the problem (\ref{L2.2}), 
since if $u_{n,\epsilon} \in W^{1,1+\epsilon}_0(\Omega)$, 
so is $|u_{n,\epsilon}|$ and

$$|\nabla |u_{n,\epsilon}||=|\nabla u_{n,\epsilon}|.$$
By regularity results, as developped by Guedda-Veron \cite{GV},
(see also Tolksdorf \cite{Tolk}), one gets that 
$u_{n,\epsilon} \in \mathcal C^{1,\alpha}(\overline{\Omega}), \alpha \in (0,1)$ 
and by Vazquez' Strict Maximum Principle \cite{Vaz}, one gets that $u_{n,\epsilon}>0$ 
in $\Omega$.
Let  $u_{n,\epsilon}$ be a solution of (\ref{L2.2}) which is positive,
then $\sigma_{n,\epsilon} =|\nabla u_{n,\epsilon}|^{\epsilon-1}
\nabla u_{n,\epsilon}$ satisfies the following partial differential equations:

\begin{equation}\label{L2.3}
\begin{cases}
-{\rm div} \sigma_{n,\epsilon} +
2n \left(\int_ \Omega u_{n,\epsilon}^{1+\epsilon} -1 \right)  
u_{n,\epsilon}^\epsilon=0 \quad \text{in} \ \Omega, \\
\sigma_{n,\epsilon}.\nabla u_{n,\epsilon} =
|\nabla u_{n,\epsilon}|^{1+\epsilon}
\quad \text{in} \ \Omega, \\
 u_{n,\epsilon} > 0, \ u_{n,\epsilon}=0 
\quad \text{\rm on} \ \partial\Omega, 
\end{cases}
\end{equation}

\begin{remark}
The solution $u_{n,\epsilon}$ is unique.\\ 
Indeed, let $u$ and $v$ two positive solutions of (\ref{L2.3}), then we have:
\begin{equation}\label{L2.4}
-{\rm div} \sigma_{\epsilon} (u)+
2n \left( \int_\Omega u^{1+\epsilon}
 -1 \right) u^{\epsilon}=0 .
\end{equation}

\begin{equation}\label{L2.5}
-{\rm div} \sigma_{\epsilon}(v)+
2n \left( \int_\Omega v^{1+\epsilon}
 -1 \right) v^{\epsilon}=0 .
\end{equation}
Let us note:
$$\alpha(u)=\int_\Omega u^{1+\epsilon} -1.$$
Substracting (\ref{L2.5}) from (\ref{L2.4}), one gets:
\begin{equation}\label{L2.6}
-{\rm div}\left(\sigma_{\epsilon}(u)-\sigma_{\epsilon}(v)\right)+
2n\left(\alpha(u)-\alpha(v)\right)u^{\epsilon}+2n \  \alpha(v)\left(u^{\epsilon}-v^{\epsilon}\right)=0 .
\end{equation}

\bigskip

\noindent{\bf Case 1}: 
$\|u\|_{1+\epsilon}^{1+\epsilon}=\|v\|_{1+\epsilon}^{1+\epsilon} .$\\
Let us multiply (\ref{L2.6}) by $(u-v)$ then integrate, we get that:

$$\int_\Omega(\sigma_{\epsilon}(u)-\sigma_{\epsilon}(v)) \cdot \nabla(u-v)+2n \ \alpha(v) \int_\Omega \left(u^{\epsilon}-v^{\epsilon}\right)(u-v)=0.$$
We know that 

\begin{equation}\label{L2.7}
\int_\Omega(\sigma_{\epsilon}(u)-\sigma_{\epsilon}(v)) \cdot \nabla(u-v)\geq 0.
\end{equation}
On the other hand it is clear that

\begin{equation}\label{L2.8} 
\int_\Omega \left(u^{\epsilon}-v^{\epsilon}\right)(u-v)\geq 0.
\end{equation}
So, we can conclude that  $\int_\Omega \left(u^{\epsilon}-v^{\epsilon}\right)(u-v)=0$,and this implies that $u=v$ a.e.

\bigskip

\noindent{\bf Case 2}:
$\|u\|_{1+\epsilon}^{1+\epsilon} \geq\|v\|_{1+\epsilon}^{1+\epsilon} .$\\   
Let us multiply  (\ref{L2.6}) by $(u-v)^+$ then integrate.
It is clear that 

$$2n\left(\alpha(u)-\alpha(v)\right)\int_\Omega u^{\epsilon}(u-v)^+ \geq 0 .$$
So we get that:

\begin{equation}\label{L2.9}
\int_\Omega(\sigma_{\epsilon}(u)-\sigma_{\epsilon}(v)) \cdot \nabla(u-v)^+
+2n \  \alpha(v)\int_\Omega\left(u^{\epsilon}-v^{\epsilon}\right)(u-v)^+ \leq0 .
\end{equation}
By (\ref{L2.7}) and  (\ref{L2.8}), we have that:
\begin{equation}\label{L2.10}
\int_\Omega(\sigma_{\epsilon}(u)-\sigma_{\epsilon}(v)) \cdot \nabla(u-v)^+
+2n \  \alpha(v)\int_\Omega \left(u^{\epsilon}-v^{\epsilon}\right)(u-v)^+ \geq0 .
\end{equation}
So from (\ref{L2.9}) and  (\ref{L2.10}), we obtain that
$$\int_\Omega(\sigma_{\epsilon}(u)-\sigma_{\epsilon}(v)) \cdot \nabla(u-v)^+
+2n \  \alpha(v)\int_\Omega\left(u^{\epsilon}-v^{\epsilon}\right)(u-v)^+ =0 .$$
Then we can conclude that 
$\int_\Omega \left(u^{\epsilon}-v^{\epsilon}\right)(u-v)^+=0$, 
and this implies that  $(u-v)^+=0$, which gives us that $u\leq v$.\\
Recall that $\|u\|_{1+\epsilon}^{1+\epsilon} \geq\|v\|_{1+\epsilon}^{1+\epsilon}$,
hence $u=v$ a.e.

\bigskip

\noindent{\bf Case 3}:
$\|u\|_{1+\epsilon}^{1+\epsilon} \leq\|v\|_{1+\epsilon}^{1+\epsilon}$.\\
We use the same arguments then the Case 2, replacing $(u-v)^+$ by $(v-u)^+$. 
  
\end{remark}

\noindent{\bf Step 2}:
\begin{proposition}\label{prop:L2}
$$\overline{\lim_{\epsilon \rightarrow 0}} 
\lambda_{n,\epsilon}(\Omega) \leq \lambda_n(\Omega).$$
\end{proposition}

\begin{proof}
Let $I_{n,\epsilon}(\phi) =\int_\Omega |\nabla \phi |^{1+\epsilon} +  
n \left(\int_\Omega| \phi|^{1+\epsilon} -1 \right)^2$, and $\delta>0$ 
be given and $\phi \in \mathcal D(\Omega)$ such that

$$\int_\Omega |\nabla \phi | +  n \left( \int_\Omega| \phi|-1 \right)^2 \leq \lambda_n+\delta,$$
For $\epsilon$ close to $0$, $|I_{n,\epsilon}(\phi)-I_{n,0}(\phi)|<\delta$, 
hence

$$\overline{\lim_{\epsilon \rightarrow 0}} 
\lambda_{n,\epsilon} \leq \lambda_n+\delta,$$
$\delta$ being arbitrary, we get 
$\overline{\lim}_{\epsilon \rightarrow 0} 
\lambda_{n,\epsilon} \leq \lambda_n.$ 
\end{proof}

\noindent Let now $u_{n,\epsilon}$ be a positive solution of (\ref{L2.3}). 
Then, it is bounded in $W^{1,1+\epsilon}_0(\Omega)$. 
Using H\"older's inequality, we get

\begin{equation}\label{L2.11}
\int_\Omega u_{n,\epsilon} dx \leq \left(\int_\Omega 
u_{n,\epsilon}^{1+\epsilon} dx\right)^{\frac{1}{1+\epsilon}} 
|\Omega|^{\frac{\epsilon}{1+\epsilon}},
\end{equation}
and thus  $(u_{n,\epsilon})_{\epsilon>0}$ is bounded in $L^1(\Omega)$. 
By the same arguments we prove that $(| \nabla u_{n,\epsilon}|)_{\epsilon>0}$ 
is also bounded in $L^1(\Omega)$. Hence, $(u_{n,\epsilon})_{\epsilon>0}$ 
is bounded in $BV(\Omega)$. 
Therefore, we may extract from it a subsequence, 
still denoted $(u_{n,\epsilon})$, 
such that

$$u_{n,\epsilon} \rightarrow u_n  \quad \text{in} \ L^k (\Omega), \forall k<1^*,$$

$$\nabla u_{n,\epsilon} \rightharpoonup  \nabla u_n 
 \quad \text{in} \ M^1(\Omega) \quad \text{weakly},$$

We need now to recall a result of concentration compactness, which is a 
consequence of the concentration compactness theory of P.L.Lions \cite{PLL}.

\begin{lemma}\label{lemma:L1}
Suppose that $\Omega$ is an open bounded set in $\R^N, N>1$, and that 
$u_{n,\epsilon}$ is bounded in $W^{1,1+\epsilon}_0(\Omega)$, then if 
$u_{n,\epsilon} \rightharpoonup u_n  \in BV(\Omega)$ weakly, there exists some 
nonnegative bounded measure $\mu$ on $\Omega$, a numerable set 
$\{x_i\}_{i \in \N} \in \overline{\Omega}$, and some numbers $\mu_i \geq 0$ such that

\begin{equation}\label{L2.12}
| \nabla u_{n,\epsilon}| \rightharpoonup \mu \geq  | \nabla u_n|+\sum_i \mu_i 
\delta_{x_i} \quad \ {\rm in} \ M^1(\Omega) \quad {\rm weakly}, 
\end{equation}
where $\delta_{x_i}$ denotes the Dirac mass on $x_i$. 
\end{lemma}

\noindent{\bf Step 3}:\\
 {\bf we obtain $\sigma_n=''\frac {\nabla u_n}{| \nabla u_n|}''$
as the weak limit of $\sigma_{n,\epsilon}=| \nabla u_{n,\epsilon}|^{\epsilon-1} \nabla u_{n,\epsilon}.$}\\
Let $\sigma_{n,\epsilon}=| \nabla u_{n,\epsilon}|^{\epsilon-1} \nabla u_{n,\epsilon}$. Then  $\sigma_{n,\epsilon}$ belongs to $L^{\frac{1+\epsilon}{\epsilon}}(\Omega)$.
By passing to the limit when $\epsilon$ goes to $0$, one obtains that $\sigma_{n,\epsilon}$ tends to $\sigma_n$ weakly in $L^{q}(\Omega)$,  $\text{for all} \  q< \infty$. We need to prove that $|\sigma_n|_{\infty} \leq 1$. 

\noindent For that aim, let $\eta$ be in $\mathcal D(\Omega, \R^N)$. Then

\begin{align*}
\left| \int_\Omega \sigma_n.\eta \right| \leq {\underline{\lim}}_{\epsilon\rightarrow 0}\left|\int_\Omega \sigma_{n,\epsilon}.\eta\right| & \leq \underline{\rm lim}_{\epsilon \rightarrow 0}\int_\Omega|\nabla u_{n,\epsilon}|^{\epsilon}|\eta|\\
&\leq \underline{\rm lim}_{\epsilon \rightarrow 0} \left(\int_\Omega| \nabla u_{n,\epsilon}|^{1+\epsilon} \right)^{\frac{\epsilon}{1+\epsilon}} \left(\int_\Omega|\eta|^{1+\epsilon} \right)^{\frac{1}{1+\epsilon}}\\
&\leq \underline{\rm lim}_{\epsilon \rightarrow 0} \left(C\right)^
{\frac{\epsilon}{1+\epsilon}} \left(\int_\Omega|\eta|^{1+\epsilon} \right)^
{\frac{1}{1+\epsilon}}\\
&\leq \int_\Omega |\eta|.
\end{align*}
On the other hand, we prove that $u_n^\epsilon$ converges weakly to some $\alpha$ in $L^\infty(\Omega),  \ \alpha \in [0, 1]$.
Indeed, for $\eta \in \mathcal D(\Omega)$,
\begin{align*}
\left|\int_\Omega u_n^\epsilon.\eta \right| 
&\leq \left(\int_\Omega|u_n^\epsilon|^{\frac{1+\epsilon}{\epsilon}} \right)^{\frac{\epsilon}{1+\epsilon}} \left(\int_\Omega|\eta|^{1+\epsilon} \right)^{\frac{1}{1+\epsilon}}\\
&\leq {\rm lim}_{\epsilon \rightarrow 0} \left(C\right)^
{\frac{\epsilon}{1+\epsilon}} \left(\int_\Omega|\eta|^{1+\epsilon} \right)^
{\frac{1}{1+\epsilon}}\\
&\leq \int_\Omega |\eta|.
\end{align*}
This implies that  $|\sigma_n|_{\infty} \leq 1$. On the other hand by passing 
to the limit in (\ref{L2.3}), one gets:

\begin{equation}\label{L2.13}
-{\rm div} \sigma_n +2n \left( \int_ \Omega u_n -1 \right)\alpha =0 .
\end{equation}

\noindent {\bf Step 4}:\\
{\bf Extension of $u_{n,\epsilon}$ outside $\Omega$ and convergence towards 
a solution of (\ref{L2.3})}.\\
We shall need in this part the Proposition \ref{prop:L1} and a classical result 
in the theory of $BV-$functions:

\begin{lemma}\label{lemma:L2}
Assume that  $\Omega$ is an open bounded set in $\R^N, \ N>1$, and that  
$u \in BV(\Omega)$. Then if $x_0 \in \Omega$, $|\sigma.\nabla u|(\{x_0\})=0$. 
\end{lemma}
\noindent The proof of this result can be found in \cite{Giu}.

\bigskip

Let $\tilde u_{n,\epsilon}$ be the extension of $u_{n,\epsilon}$ by $0$ in 
$\R^N- \overline{\Omega}$. 
Then $\tilde u_{n,\epsilon} \in W^{1,1+\epsilon}(\R^N)$, since $u_{n,\epsilon}=0$ on $\partial \Omega$, and $(\tilde u_{n,\epsilon})$ is bounded in  $BV(\R^N)$.
Then one may extract from it a subsequence, still denoted $(\tilde u_{n,\epsilon})$ such that

$$\tilde u_{n,\epsilon} \rightarrow v_n  \quad \text{in} \ L^k (\R^N), \forall k< \frac{N}{N-1} , $$
with $v_n=0$ outside of $\overline{\Omega}$. We denote by $u_n$ the restriction of $v_n$ to  $\Omega$. In addition: 

$$\nabla \tilde u_{n,\epsilon} \rightharpoonup  \nabla v_n  \quad \text{in} \ M^1(\R^N) \quad \text{weakly},$$

$$\sigma_{n,\epsilon}=|\nabla u_{n,\epsilon}|^{\epsilon-1} \nabla u_{n,\epsilon} \rightharpoonup \sigma_n \quad \text{in} \ L^q(\Omega), \forall q < \infty . $$
Using Concentration Compactness Lemma, there exists a non negative measure $\mu$, with support in $\overline{\Omega}$, a numerable set $\{x_i\}_i$ in $\overline{\Omega}$ and some non negative reals $\mu_i$, such that 

$$|\nabla \tilde u_{n,\epsilon}|^{1+\epsilon} \rightharpoonup \mu \geq |\nabla v_n|+ {\Sigma}_i \mu_i\delta_{x_i} . $$
Multiplying (\ref{L2.3}) by $\tilde u_{n,\epsilon} \phi$ where $\phi \in \ \mathcal D(\R^N)$, and integrating by parts, one obtains:

$$\int_{\overline{\Omega}} \sigma_{n,\epsilon} \cdot \nabla (\tilde u_{n,\epsilon} \phi)+2n \left( \int_ \Omega \tilde u_{n,\epsilon}^{1+\epsilon} -1\right)\int_ \Omega \tilde u_{n,\epsilon}^{1+\epsilon} \phi =0 , $$
or equivalently 

$$\int_{\R^N} |\nabla (\tilde u_{n,\epsilon})|^{1+\epsilon} \phi+\int_{\R^N}\sigma_{n,\epsilon} \tilde u_{n,\epsilon} \cdot \nabla \phi+2n \left(\int_{\R^N}\tilde u_{n,\epsilon}^{1+\epsilon} -1\right) \int_{\R^N}\tilde u_{n,\epsilon}^{1+\epsilon} \phi =0 .$$
Since $\sigma_{n,\epsilon} \rightharpoonup \sigma_n$ in  $L^q (\Omega)$ for all $q$,
in particular for some $\alpha >0$, $\sigma_{n,\epsilon}$ tends weakly towards $\sigma_n$ in  $L^{N+\alpha} (\Omega)$ , and then, since $\tilde u_{n,\epsilon}$ tends strongly towards $v_n$ in $L^k, k< \frac{N}{N-1}$, one obtains that:

$$\int_{\R^N}\sigma_{n,\epsilon} \tilde u_{n,\epsilon} \cdot \nabla \phi \longrightarrow \int_{\R^N}\sigma_n v_n \cdot \nabla \phi,  \ \quad \text{when} \ \epsilon\rightarrow 0 . $$
By passing to the limit in the last equation above, one obtains:

\begin{equation}\label{L2.14}
\langle \mu, \phi \rangle+ \int_{\Omega}\sigma_n u_n \cdot \nabla \phi+2n \left(\int_{\Omega}u_n-1\right)\int_{\Omega}u_n \phi =0 .
\end{equation}
Using generalised Green's Formula in Proposition \ref{prop:L1} and (\ref{L2.13}), we have

\begin{align}\label{L2.15}
\int_{\Omega}\sigma_n u_n \cdot \nabla \phi=
&-\int_{\Omega} {\rm div} \sigma_n u_n \phi-\int_{\Omega}\sigma_n \cdot \nabla u_n \phi+\int_{\partial \Omega} \sigma_n \cdot \overrightarrow {n} u_n \ \phi \nonumber\\
&=-2n \left(\int_{\Omega}u_n -1\right)\int_{\Omega} u_n \phi- \int_{\overline{\Omega}}\sigma_n \cdot \nabla v_n  \phi.
\end{align}
Substracting (\ref{L2.15}) from (\ref{L2.14}), one gets for $\phi \in \mathcal D(\overline{\Omega})$\begin{equation}\label{L2.16}
\langle \mu, \phi \rangle- \int_{\overline{\Omega}}\sigma_n \cdot \nabla v_n \phi=0 .
\end{equation}
Let now $h$ be a $|\nabla v_n|-$measurable function and $\mu^\bot$ be a measure orthogonal to $|\nabla u_n|$, such that, according to the Radon-Nikodym decomposition, one has
 
\begin{equation}\label{L2.17}
\mu= h |\nabla v_n| +\mu^\bot .
\end{equation}
By Lemma \ref{lemma:L2} and the analogous of (\ref{L2.12}) of Lemma \ref{lemma:L1}, one has

\begin{equation}\label{L2.18} 
h |\nabla v_n| \geq |\nabla v_n|,
\end{equation}
and 

\begin{equation}\label{L2.19}
\mu^\bot \geq \sum_i \mu_i\delta_{x_i} .
\end{equation}
Using (\ref{L2.18}) and (\ref{L2.19}) in equation (\ref{L2.16}) one gets that

\begin{equation}\label{L2.20}
h |\nabla v_n|=\sigma_n \cdot \nabla v_n   \ \quad \text{on} \ \overline{\Omega},
\end{equation}
and

\begin{equation}\label{L2.21}
\mu^\bot= \sum_i \mu_i\delta_{x_i} \leq \sigma_n \cdot \nabla v_n .
\end{equation}
Using Lemma \ref{lemma:L2} one gets that $\mu_i=0 \ \forall i$. 
And from (\ref{L2.20}), $|\sigma_n| \leq 1$ and Proposition \ref{prop:L1}, one obtains that in the sense of measures:
$$|\sigma_n \cdot \nabla v_n| \leq |\nabla v_n|  \ \quad \text{on} \ \overline{\Omega} , $$

$$\mu=\sigma_n \cdot \nabla v_n , $$ 
and then 

$$\sigma_n \cdot \nabla v_n=|\nabla v_n| \ \quad \text{on} \ \overline{\Omega} , $$
and $h=1$, $|\nabla v_n|-$almost everywhere. 
Using this in equation (\ref{L2.14}) with $\phi=1$, we have:

$$\int_{\Omega}|\nabla v_n|+\int_{\partial \Omega} v_n+2n \left(\int_{\Omega}v_n -1\right)\int_{\Omega}v_n=0 . $$
Recalling that we have from Proposition \ref{prop:L1} that:

$$\nabla v_n = \nabla v_n \chi_{\Omega} -v_n  \ \delta_{\partial \Omega} \overrightarrow {n} , $$

$$\sigma_n \cdot \nabla v_n=\sigma_n \cdot \nabla v_n \chi_{\Omega} -\sigma_n.\overrightarrow {n} v_n \delta_{\partial \Omega} . $$
This implies that $\sigma_n.\nabla v_n=|\nabla v_n|$, on $\Omega \cup \partial \Omega$. This condition can be splitted in the two equations

$$\begin{cases}
\sigma_n \cdot \nabla u_n=|\nabla u_n|  \ \quad \text{in} \  \Omega , \\
\sigma_n \cdot \overrightarrow {n} u_n = - u_n \ \quad \text{on} \ \partial \Omega . 
\end{cases}$$
Then $u_n$ is a nonnegative solution of (\ref{L1.2}). Moreover, the convergence of $|\nabla \tilde u_{n,\epsilon}|$ is tight on $\overline{\Omega}$ which means that 

$$\int_{\Omega}|\nabla  u_{n,\epsilon}| \longrightarrow\int_{\Omega}|\nabla u_n|+\int_{\partial \Omega} u_n, \ \quad \text{when} \ \epsilon\rightarrow 0 . $$
Indeed, one has $\int_\Omega |\nabla u_{n,\epsilon}|^{1+\epsilon}\rightarrow \int_\Omega |\nabla u_n|+\int_{\partial \Omega }u_n$ and using the lower semicontinuity for the extension $v_{n,\epsilon}$,  we get first 

$$ \int_\Omega |\nabla u_n|+\int_{\partial \Omega }u_n\leq \underline{\lim}_{\epsilon\rightarrow 0} \int_{\Omega} |\nabla u_{n,\epsilon}|,$$
and secondly, using  H\"older's inequality $\int_\Omega |\nabla u_{n,\epsilon}|^{1+\epsilon}\geq (\int_\Omega|\nabla u_{n,\epsilon}|) ^{\frac{1}{1+\epsilon}} |\Omega |^{\frac{-\epsilon}{1+\epsilon}}$ which gives by passing to the limit the reverse inequality.

\noindent{\bf Step 5}:
{\bf $u_n$ is a solution of (\ref{L1.4})}\\
Let us recall the relaxed form of (\ref{L1.3})

$$\inf_{\begin{smallmatrix} 
u \in BV(\Omega) 
\end{smallmatrix}} 
\quad \left\{\int_\Omega |\nabla u | +\int_{\partial \Omega} |u|+
n \left( \int_\Omega|u| -1 \right)^2 \right \} . \ \quad \  \eqno{(1.4)}$$
We prove now that the solution $u_n$ obtained in the two previous steps 
is a nonnegative solution of (\ref{L1.4}). 
For that aim, let us recall that by using the lower semi-continuity
for the weak topology of $BV(R^N)$, we have:

\begin{align*}
\lambda_n(\Omega)=&\int_{\R^n}|\nabla v_n|+n\left(\int_{\Omega} v_n -1 \right)^2 \\
&\leq \underline{\rm lim}_{\epsilon \rightarrow 0}\int_{\R^n}
|\nabla \tilde u_{n,\epsilon}|^{1+\epsilon}
+n\left(\int_{\Omega}\tilde u_{n,\epsilon}^{1+\epsilon} -1\right)^2 \\ 
&\leq \overline{\rm lim}_{\epsilon \rightarrow 0} \lambda_{n,\epsilon}(\Omega)
\leq \lambda_n(\Omega) .
\end{align*}
Using the fact  that 
$$\int_{\R^n}|\nabla v_n|=\int_{\Omega \cup \partial \Omega}|\nabla v_n|=
\int_\Omega |\nabla u_n|+\int_{\partial \Omega} u_n . $$
One obtains that $u_n$ is a nonnegative solution of the relaxed problem (\ref{L1.4})
and in the same time we get

$$\lim_{\epsilon \rightarrow 0} \lambda_{n,\epsilon} = \lambda_n . $$ 
Moreover, one has:

$$\int_{\Omega} |u_{n,\epsilon}| \longrightarrow \int_{\Omega}  u_n,
\quad \text{when} \ \epsilon \rightarrow 0 . $$
Then, we conclude that:

$$\int_{\Omega} |\nabla u_{n,\epsilon}|\longrightarrow 
\int_{\Omega} |\nabla u_n|+\int_{\partial \Omega} u_n,
\quad \text{when} \ \epsilon \rightarrow 0 . $$
Hence we get the tight convergence of $u_{n,\epsilon}$ towards $u_n$ in 
$BV(\overline{\Omega})$.

\end{proof}

\section{Convergence result}

We begin to recall some properties of the first eingenvalue for the 1-Laplacian
operator, see e.g. \cite{DF3}

\begin{proposition}\label{prop:L3}
Suppose that $\lambda > 0$ is such that there exists $\sigma, \ |\sigma|_{\infty} \leq 1$, and $u \geq 0$ in $BV(\Omega)$ with

\begin{equation}\label{L3.1}
\begin{cases}
-{\rm div} \sigma=\lambda, \ u \geq0, \ u \not \equiv 0, \ u \in BV(\Omega),\\
\sigma.\nabla u= |\nabla u|  \quad \text{\rm in} \ \Omega, \ |\sigma|_{L^{\infty}( \Omega)} \leq 1, \\
\sigma.\overrightarrow {n} (-u)=u \ \quad \text{\rm on}\ \partial\Omega .
\end{cases}
\end{equation}
Then $\lambda=\lambda_1$ where

$$\lambda_1=\inf_{\begin{smallmatrix} u \in W^{1,1}_0(\Omega) \\ \| u\|_1 =1\end{smallmatrix}} \quad \int_\Omega |\nabla u| , $$
$\lambda_1$ is called the first eingenvalue for $-\Delta_1$ on $\Omega$. 
Moreover

\begin{equation}\label{L3.2}
\lambda_1=\inf_{\begin{smallmatrix} u \in BV(\Omega) \\ 
\| u\|_1 =1\end{smallmatrix}} 
\quad \left\{\int_\Omega |\nabla u | +\int_{\partial \Omega} |u|\right\} ,
\end{equation}
and this last infimum is achieved on some $u$ which satisfies {\rm (\ref{L3.1})}. Among the ``eingenfunctions'' there exist caracteristic functions of Cacciopoli sets.
\end{proposition}

\begin{theorem}\label{th:L2}
Let $u_n$ be a function for which $\lambda_n$ is achieved, then, up to a subsequence, $(u_n)$ converges to $u \in BV(\Omega), \ u \geq0, \ u \not \equiv 0,$ which realizes the minimum defined in {\rm (\ref{L3.2})}.
Moreover

$$\lim_{n \rightarrow \infty} \lambda_n = \lambda_1 . $$ 
\end{theorem}

\begin{proof}[Proof of the Theorem \ref{th:L2}] 
For $\lambda_n$ and $\lambda_1$ defined as above, it is clear that we have:

\begin{equation}\label{L3.3}
\overline{\rm lim}_{n \rightarrow \infty} \lambda_n \leq \lambda_1.
\end{equation}
Let $(u_n)_n$ be a sequence of solutions  of the relaxed problem (\ref{L1.4}). We begin to prove that $(u_n)_n$ is bounded in $BV(\Omega)$.
For that aim let us note that by (\ref{L3.3}), one gets that  
$n\left(\int_\Omega u_n -1 \right)^2$ is bounded by $\lambda_1$, which implies that ${\rm lim}_{n \rightarrow \infty} \left(\int_\Omega u_n-1\right)^2 = 0$.
Then

$${\rm lim}_{n \rightarrow \infty} \int_\Omega u_n = 1 , $$
Hence, $(u_n)_n$ is bounded in $L^1(\Omega)$.\\
Using once more (\ref{L3.3}) for $|\nabla(u_n)|_n$ one can conclude that $(u_n)_n$ is bounded in $BV(\Omega)$.
\noindent Then, the extension of $u_n$ by zero outside of $\overline{\Omega}$ is bounded in $BV(\R^N)$. One can then extract from it a subsequence, still denoted $u_n$, such that 
$$u_n \rightharpoonup u  \quad \text{in} \ BV(\R^N)\quad \text{weakly},$$
Obviously $u=0$ outside of $\overline{\Omega}$. By compactness of the Sobolev embedding from $BV(\Omega)$ into $L^1(\Omega)$, one has $|u|_{L^1(\Omega)}=1$. Using lower semi continuity, one has

\begin{align*}
\lambda_1 \leq \int_{\R^N} |\nabla u| &\leq  \int_{\R^N}|\nabla u |+ n\left(\int_{\R^N}u -1\right)^2 \\
&\leq \underline{\rm lim}_{n \rightarrow \infty} \left(\int_{\R^N}|\nabla u_n|+ n\left(\int_{\R^N} u_n -1\right)^2 \right) \\
&\leq \overline{\rm lim}_{n \rightarrow \infty} \lambda_n \leq \lambda_1.
\end{align*}
Then one gets that 

$${\rm lim}_{n \rightarrow \infty} \lambda_n =\lambda_1 . $$
Since $u=0$ outside of $\overline{\Omega}$, one has $\nabla u =\nabla u \chi_{\Omega} - u \cdot \overrightarrow {n} \ \delta_{\partial \Omega}$ on $\partial \Omega$,
and then

$$\int_{\R^N}|\nabla u |=\int_{\Omega}|\nabla u|+\int_{\partial \Omega}u . $$
Moreover, one obtains that:

$${\rm lim}_{n \rightarrow \infty}{\rm n}\left( \int_\Omega u_n-1\right)^2=0 , $$
and 

$${\rm lim}_{n \rightarrow \infty} \int_\Omega|\nabla u_n|= \int_\Omega|\nabla u| . $$
Then, we get the tight convergence of $u_n$ to $u$ in $BV(\overline{\Omega})$.
 
Hence by passing to the limit in (\ref{L2.13}) when $n \rightarrow \infty$ one gets that:
$$-2n \left( \int_ \Omega u_n -1 \right) \longrightarrow \lambda_1, \ \quad \text{when} \ \ n \rightarrow \infty.$$
\end{proof}

\section{Minima as Cacciopoli sets}

Let us introduce $\lambda_{0,n}$ as the value of the infimum

\begin{equation}\label{L4.1}
\lambda_{0,n}=\inf_{\begin{smallmatrix}E,E \ \text{is Cacciopoli set}\\E \subset \subset \Omega 
\end{smallmatrix}} 
\left\{ \int_ \Omega |\nabla \chi _E| + n\left( \int_ \Omega\chi _E-1 \right)^2 \right\},
\end{equation}
(let us recall that a Cacciopoli set in $\Omega$ is merely a set whose caracteristicfunction belongs to $BV(\Omega)$).
We have the following.

\begin{theorem}\label{th:L3}
One has
\begin{align*}\lambda_{0,n}&=\inf_{\begin{smallmatrix}E,E \  \text{is Cacciopoli set}\\E \subset \subset \Omega \end{smallmatrix}} 
\left\{ \int_ \Omega |\nabla \chi _E| + n\left( \int_ \Omega\chi _E-1 \right)^2 \right\}\\&=\inf_{\begin{smallmatrix}E,E \ \text{is Cacciopoli set}\\E \text{set in} \ \R^n \end{smallmatrix}} \left\{ \int_ \Omega |\nabla \chi _E| +\int_ {\partial \Omega} \chi _E+ n\left( \int_ \Omega\chi _E-1 \right)^2 \right\}\\&=\inf_{\begin{smallmatrix}E,E \ \text{is Cacciopoli set}\\E \text{set in} \ \R^n \end{smallmatrix}}P(E,\Omega)+|E \cap \partial \Omega|+n\left(|E \cap \Omega|- 1\right)^2\end{align*}and
$$\lambda_{0,n} \geq \lambda_{1,n} . $$
\end{theorem}

\begin{remark}
$P(E,\Omega)$ is the perimeter of E in $\Omega$, see (\cite{DGCL}, \cite{GMS}).
\end{remark}

\begin{proposition}\label{prop:L4}
Let $\lambda_{0,n}$ be defined as in Theorem \ref{th:L3}, then $\lambda_{0,n}$ is achieved.
\end{proposition}

\begin{proof} 
Let $(E_i)$ be a subsequence of Cacciopoli sets, $E_i\subset \subset \Omega$ such that
$$\int_ \Omega |\nabla \chi _{E_i}| + n\left( \int_ \Omega\chi _{E_i}-1 \right)^2 \rightarrow \lambda_{0,n} . $$
Then
$$\overline \lim_{i\rightarrow \infty} \left\{ \int_\Omega |\nabla \chi _{E_i}| + n\left( \int_ \Omega\chi _{E_i}-1 \right)^2\right\} \leq \lambda_{0,n} . $$
It is clear that $\chi _{E_i}$ is bounded in $BV(\Omega)$ (same arguments as in Theorem \ref{th:L3}).
More precisally $\chi _{E_i}$ is bounded in $BV(\R^n)$.\\Extracting from it a subsequence still denoted $\chi _{E_i}$, one get that
$$\chi _{E_i} \rightharpoonup u \ \text{in}BV(\R^n) . $$
By construction $u=0$ outside of $\overline{\Omega}$. Moreover one can assume that$\chi _{E_i}$ tends to $u$ a.e, and then $u$ can only takes the values $0$ and $1$.
As a consequence $u$ is the caracteristic function of some set E.\\By lower semicontinuity, one has that 
$$\int_{\R^n} |\nabla \chi _{E}| + n\left( \int_{\R^n} \chi _{E}-1 \right)^2 \leq \underline \lim_{i\rightarrow \infty} \left\{ \int_\Omega |\nabla \chi _{E_i}| + n\left( \int_ \Omega\chi _{E_i}-1 \right)^2\right\} , $$
Then, one obtains that $E$ is a solution for the relaxed problem (\ref{L4.1}). \end{proof}


\addcontentsline{toc}{section}{Bibliographie}
\begin{thebibliography}{ab}
 
\bibitem{ACCh1}{\bf Alter, F. ; Cazelles, V. ; Chambolle, A.} {\em A characterization of convex calibrable sets in $\R^N$}, prepublication.
 
\bibitem{ACCh2}{\bf Alter, F. ; Cazelles, V. ; Chambolle, A.} {\em Evolution of convex sets in the plane by the minimizing totalvariaion flow }, Prépublication.  
 
\bibitem{ACM} {\bf Andreu, F. ; Caselles, V. ; Mazón, J. M.} {\em A strongly degenerate quasilinear elliptic equation},  Nonlinear Anal.  61  (2005),  no. 4, 637--669
 
\bibitem{BCN}{\bf Bellettini, G. ; Caselles, V. ; Novaga, M.} {\em Explicit solutions of the eigenvalue problem $-{\rm div} (\frac{Du}{|Du|}) = u$.} 
 
\bibitem {Ch} {\bf Cheeger, J.} {\em A lower bound for the smallesteigenvalue of the Laplacian} in Problems in Analysis, Symposium in honorof Salomon Bochner, Ed : RC Ganning, Princeton Univ. Press (1970), p.195-199. 

\bibitem{DGCL} {\bf De Giorgi; Carriero, M.; Leaci, A.;} {\em Existence Theorem for a minimum problem with a Free dicountinuity set,}A.R.M.A, 108, 195-218, (1989).

\bibitem{DF1} {\bf Demengel, F.} {\em On Some Nonlinear Partial Differential Equations Involving The 1-Laplacian and Critical Sobolev exponent},  ESAIM: Control, Optimisation and Calculus of Variations, 4 (1999), 667-686
 
\bibitem{DF2} {\bf Demengel, F.} {\em Some compactness result for some spaces of functions with bounded derivatives,} \ A.R.M.A. 105(2) (1989), 123-161.
  
\bibitem{DF3} {\bf Demengel, F.} {\em Theoremes d'existence pour des equations avec l'operateur $1$-Laplacien, premiere valeur propore pour -$\Delta_1$,} \ C.R Acad. Sci. Paris, Ser. I334 (2002), 1071-1076.

\bibitem{DF4} {\bf Demengel, F.} {\em Some existence's results for noncoercive $1-$Laplacian operator,} Asymptot. Anal.  43  (2005),  no. 4, 287-322.

\bibitem{DF5} {\bf Demengel, F.} {\em Functions locally almost $1-$harmonic},
Applicable Analysis, Vol.83, N°9, September 2004, 865-896.

\bibitem{DF6} {\bf Demengel, F.}  {\em On some nonlinear partial differential equations involving the "$1-$Laplacian " and critical Sobolev exponent,}  ESAIM Control Optim. Calc. Var.  4  (1999), 667--686.

\bibitem{Ek.Tem}  {\bf Ekeland, I. and Temam, R.} {\em Convex Analysis and variational problems,}
North-Holland, 1976.

\bibitem{Giu}  {\bf Giusti, E.} {\em Minimal surfaces and functions of bounded variation,}Notes de cours rediges pr G.H. Williams.Departement of Mathematics Australian National University, Canberra (1977), et Birkhauser (1984).
 
\bibitem{GMS} {\bf Giaquinta, M. ; Modica, G. ; and Soucek, J.} {\em Cartesian Currents in the Calculus of Variations I,}LNM, Vol 37, Springer, 1997.
 
\bibitem{GV} {\bf Guedda, M. and Veron, L.} {\em Quasilinear elliptic equations involving critical sobolev exponents,} Nonlinear Analysis, Theory, Methods and Applications, 13 (1989), 879-902.

\bibitem{Ko.Tem}  {\bf Kohn, R.V. and Temam, R.} {\em Dual spaces of stress and strains with applications to Hencky plasticity,}Appl.Math.Optim (10) (1983), 1-35.

\bibitem{PLL} {\bf Lions, P.L.} {\em  The concentration-compactness principle in the calculus of variations. The limit case,} I et II.  Rev. Mat. Iberoamericana  1  (1985),  no. 1, 145--201. 

\bibitem{Str.Tem}  {\bf Strang, G. and Temam, R.} {\em Functions with bounded derivatives,} A.R.M.A. (1980), 493-527.

\bibitem{Tolk} {\bf Tolksdorf, P.} {\em Regularity for a more general class of quasilinear elliptic equations,} Journal of Differential Equations, 51 (1984), 126-150.
 
\bibitem{Vaz} {\bf  Vazquez, J.L.} {\em A Strong maximum principle for some quasilinear elliptic equations,} Appl. Math. Optim., 12 (1984), 191-202.
 
 
\end{thebibliography}
\end{document}